\documentclass[twocolumn]{autart}
\usepackage[noadjust]{cite}
\usepackage[dvips]{graphicx}
\usepackage[cmex10]{amsmath}
\usepackage{amsfonts}
\usepackage{amssymb}

\usepackage{tikz}
\usetikzlibrary{arrows, decorations.pathmorphing, backgrounds, positioning, fit, calc}

\newcommand{\bseq}{\begin{subequations}}
\newcommand{\eseq}{\end{subequations}}
\newcommand{\bmat}{\begin{bmatrix}}
\newcommand{\emat}{\end{bmatrix}}
\newcommand{\beq}{\begin{equation}}
\newcommand{\eeq}{\end{equation}}
\newcommand{\beqs}{\begin{equation*}}
\newcommand{\eeqs}{\end{equation*}}
\newcommand{\bali}{\begin{aligned}}
\newcommand{\eali}{\end{aligned}}
\newcommand{\RR}{\mathbb{R}}
\newcommand{\CC}{\mathbb{C}}

\newcommand{\barr}{\begin{array}}
\newcommand{\earr}{\end{array}}

\newcommand{\Diag}{\textrm{Diag}}

\usepackage[position=t]{subfig}

\usepackage{epstopdf}
\usepackage{dashrule}
\newcommand{\aaand}{\qquad \mathrm{and} \qquad}

\newcommand{\iiff}{\qquad \Leftrightarrow \qquad}

\newcommand{\norm}[1]{\left\lVert#1\right\rVert}

\newtheorem{ass}{\bf{Assumption}}

\newtheorem{myeg}{\bf{Example}}

\begin{document}

\begin{frontmatter}

\title{Compressed Sensing for Network Reconstruction\thanksref{footnoteinfo}} 

\thanks[footnoteinfo]{This paper was not presented at any IFAC 
meeting. Corresponding author D.~Hayden.}
% Tel. +XXXIX-VI-mmmxxi. Fax +XXXIX-VI-mmmxxv.}

%\author[Paestum]{Marcus Tullius Cicero}\ead{cicero@senate.ir},    % Add the 
%\author[Rome, Baiae]{Julius Caesar}\ead{julius@caesar.ir},               % e-mail address 
%\author[Baiae]{Publius Maro Vergilius}\ead{vergilius@culture.ir}  % (ead) as shown
%
%\address[Paestum]{Buckingham Palace, Paestum}  % Please supply                                              
%\address[Rome]{Senate House, Rome}             % full addresses
%\address[Baiae]{The White House, Baiae}        % here.

\author[cam]{David Hayden}\ead{dph34@cam.ac.uk},
\author[ucb]{Young Hwan Chang}\ead{yhchang@berkeley.edu},
\author[cam,lux]{Jorge Goncalves}\ead{jmg77@cam.ac.uk},
\author[ucb]{Claire Tomlin}\ead{tomlin@eecs.berkeley.edu}

\address[cam]{Department of Engineering, University of Cambridge, CB2 1PZ, UK}
\address[ucb]{Department of Electrical Engineering and Computer Sciences, University of California, Berkeley, CA 94720 USA}
\address[lux]{University of Luxembourg, Facult\'{e} des Sciences, de la Technologie et de la Communication, 7 Avenue des Hauts Fourneaux, L-4362 BELVAL, Luxembourg}

% Experiment design; model selection; networked control systems; modelling and identification; dynamics and control
%Bio control; closed-loop identification; directed graphs; estimation algorithms; identifiability; interconnections matrices; linear equations; multivariable systems; network topologies; transfer function matrices.
%
%\author{David~Hayden, Young~Hwan~Chang, Jorge~Gon\c calves and Claire~Tomlin%
%\thanks{D. Hayden is with Department of Engineering, University of Cambridge, CB2 1PZ, UK, email: {dph34@cam.ac.uk}.}%
%\thanks{J. Gon\c{c}alves is with Department of Engineering, University of Cambridge, CB2 1PZ, UK and University of Luxembourg, Facult\'{e} des Sciences, de la Technologie et de la Communication, 7 Avenue des Hauts Fourneaux, L-4362 BELVAL, Luxembourg, email: {jmg77@cam.ac.uk}.}%
%\thanks{Y. H. Chang and C. Tomlin are with Department of Electrical Engineering and Computer Sciences, University of California, Berkeley, CA 94720 USA, email: {yhchang@berkeley.edu, tomlin@eecs.berkeley.edu}.}%
%\thanks{This research was supported by the Engineering and Physical Sciences Research Council under Grant EP/G066477/1 and by the NIH NCI under the ICBP and PS-OC programs (5U54CA112970-08).}}

\begin{abstract}
The problem of identifying sparse solutions for the link structure and dynamics of an unknown linear, time-invariant network is posed as finding sparse solutions $x$ to $Ax=b$. If the sensing matrix $A$ satisfies a rank condition, this problem has a unique, sparse solution. Here each row of $A$ comprises one experiment consisting of input/output measurements and cannot be freely chosen. We show that if experiments are poorly designed, the rank condition may never be satisfied, resulting in multiple solutions. We discuss experimental strategies for designing experiments such that the sensing matrix has the desired properties and the problem is therefore well posed. This formulation allows prior knowledge to be taken into account in the form of known nonzero entries of $x$, requiring fewer experiments to be performed. A number of simulated examples are given to illustrate the approach, which provides a useful strategy commensurate with the type of experiments and measurements available to biologists. We also confirm suggested limitations on the use of convex relaxations for the efficient solution of this problem.
\end{abstract}

\end{frontmatter}

%\maketitle

%%%%%%%%%%%%%%%%%%%%%%%%%%%%%%%%%%%%%%%%%%%%%%%%%%%%

\section{Introduction}

Compressed Sensing (CS) refers to the ability to find a sparse solution $x$ to the underdetermined set of equations $Ax=b$ \cite{donoho}. This problem is relevant in applications in computer vision and signal processing, where a signal often has a representation that is sparse in some domain and can hence be recovered by making relatively few samples in that domain \cite{Candes2008}. Specifically, suppose some signal $\theta \in \RR^n$ can be expressed in an orthonormal basis $\Phi \in \RR^{n \times n}$ such that $x = \Phi \theta$ where $x$ is sparse in the sense that $\|x\|_0 = k < n$ and $\| \cdot \|_0$ denotes the number of nonzero entries of a vector. By taking $m<n$ samples of $x$ via an appropriately chosen sensing matrix $A \in \RR^{m \times n}$, we may recover $x$ and hence $\theta$.

A related problem is that of identifying the link structure and dynamics of an unknown network from certain observations of it. This is a general inverse problem, currently of particular importance in cell-biological applications, such as identifying Genetic Regulatory Networks (GRNs) \cite{howtoinfer, PNAS_Review, Nature_Review, crowdwisdom, michail}. In this context the problem is typically underdetermined due to both a paucity of data and limitations on the number of experiments that can be performed. The underlying network is often known to be sparse in the sense that the degree of each node is bounded, and the assumption of sparsity is commonly used as an heuristic to obtain a solution \cite{yeung, NIR, TSNI_a, napol, antonis, sanad, Chang_CDC11, bolstad, solo, chiuso, mater}. This problem is fundamentally different from typical CS applications in that the sensing matrix cannot be chosen freely, but arises based on the experiments applied and the underlying system itself.

The problem was considered for Linear, Time-Invariant (LTI) systems with full state measurement in \cite{Chang_CDC11} using time-series data of a single perturbation to the network. If sufficiently many time points are observed, the solution is shown to be unique and hence only one such experiment is required. By assuming that the solution is sparse, the required number of time points can be reduced. In \cite{sanad}, sparse networks of FIR filters are treated, again with the aim of reducing the number of data points required. The parameters of the filters are estimated using Block Orthogonal Matching Pursuit \cite{eldar} and the notion of \emph{network coherence} is introduced, referring to the coherence (see \cite{Candes2007}) of a network-derived sensing matrix. The network coherence is observed to have a lower bound for some cases, which suggests a possible limitation to the use of CS for network reconstruction.

Other examples from the literature include \cite{chiuso}, in which the problem is posed as sparse input selection for MISO LTI systems and \cite{mater} which concerns the estimation of sparse MISO Wiener filters. For general MIMO LTI systems with deterministic inputs, it was shown in \cite{Gon08} that a certain number of targeted inputs are required in order for the problem to be well posed. This is equivalent to performing experiments to probe the network, for example in a biological context using genetic mutations to identify GRNs \cite{ECC13}. Here we suppose that sufficient data points are available and consider whether the assumption of sparsity can be used to reduce the number of such experiments required. Our focus is therefore on the \emph{identifiability} of the network, rather than a particular method, although the approach naturally provides an algorithm for steady-state or frequency-domain identification.

Our contributions are as follows: first we show that for a poor choice of experiments even the sparsest solution may not be unique; then the problem of experiment design to ensure solution uniqueness is addressed; finally, simulated examples demonstrate the effectiveness of the experiment design but reveal similar limitations to those observed in \cite{sanad} with regards to the network coherence. Hence although the sparse solution may be unique, it may be difficult to find using, for example, basis pursuit.

In Section \ref{sec:back} we review some standard results in CS and network reconstruction for LTI systems. Then in Section \ref{sec:prior} we discuss how prior knowledge can be incorporated directly into the CS framework and how this reduces the number of experiments needed for exact reconstruction. Section \ref{sec:main} addresses the problem of underdetermined network reconstruction, first showing that the standard assumptions of CS are not sufficient for exact reconstruction, then proposing experimental procedures to ensure exact reconstruction. In Section \ref{sec:sim} a number of simulation examples are presented to support the results and conclusions are given in Section \ref{sec:conc}.

\subsection*{Notation}
Denote by $A(i,j)$, $A(i,:)$ and $A(:,j)$ entry $(i,j)$, row $i$ and column $j$ respectively of matrix $A$ and by $A^T$ its transpose. The diagonal matrix comprising the diagonal entries of $A$ is denoted $\textrm{Diag}(A)$. The function $\|x\|_0$ (the $l_0$ ``norm'') returns the number of nonzero entries in the vector $x$.

\section{Background} \label{sec:back}

\subsection{Compressed Sensing}
Sparse solutions, $x \in \mathbb{R}^n$, are sought to the following problem:
\begin{equation} \label{Axb}
	Ax = b
\end{equation}

\noindent where $A \in \mathbb{R}^{m\times n}$ and $b \in \mathbb{R}^m$ are known and $m<n$. The sparsest such solution (or set of solutions) are the minimizing argument(s) of:
\beq
	\min_{x} \norm{x}_0 \qquad \textrm{subject to} \qquad Ax = b
	\label{eq:sol_CS}
\eeq

\noindent The following well-known lemma provides a sufficient condition that the solution to (\ref{eq:sol_CS}) is unique:
\begin{lem} \label{lem:cs}
If the sparsest solution to (\ref{eq:sol_CS}) has $\norm{x}_0 = k$  and $m \geq 2k$ and all subsets of $2k$ columns of $A$ are full rank, then this solution is unique.
\end{lem}

\begin{pf}
Suppose two solutions exist: $Ax^{(1)} = b$ and $Ax^{(2)} = b$, where $\|x^{(1)}\|_0 = \|x^{(2)}\|_0 = k$ and subtract one equation from the other: $A(x^{(1)} - x^{(2)}) = 0$. Since $\|x^{(1)} - x^{(2)}\|_0 \leq 2k$, this equation is equivalent to: $\hat{A}\hat{x} = 0$, where $\hat{A} \in \mathbb{R}^{m \times l}$, $0 \neq \hat{x} \in \mathbb{R}^l$, $l \leq m$ and $\hat{A}$ is full rank, which is a contradiction.
\end{pf}

Convex relaxations of \eqref{eq:sol_CS} are typically sought, such as $l_1$ minimization (basis pursuit \cite{Candes2008}), which can be solved by linear programming:
\beq
	\min_{x} \norm{x}_{1}  \qquad \textrm{subject to} \qquad Ax = b
	\label{eq:sol_CS1}
\eeq

\noindent It has been shown that $l_1$ minimization solves exactly \eqref{eq:sol_CS} if $m$ is sufficiently large and the matrix $A$ is sufficiently \emph{incoherent} \cite{donoho2, Candes2007}. The coherence of a matrix is defined as follows:
\beqs
\mu(A) = \max_{i< j} \frac{|A_i^T A_j|}{\|A_i\|_2\|A_j\|_2}
\eeqs

\noindent where $A_i$ denotes the $i^{th}$ column of $A$. Numerical simulations suggest that in practice, most $k$-sparse signals require $m\ge 4k$ in order to be recovered exactly \cite{Candes2007}.

\subsection{Dynamical Networks}

Consider a vector of directly observed variables $y(t) \in \RR^p$, whose entries $y_i(t)$ are governed by the following set of LTI equations for $i=1, \ldots p$:
\beq \label{yqp}
y_i(t) = \sum_{j \neq i} q_{ij}(t) * y_j(t) + \sum_{k=1}^r p_{ik}(t) * u_k(t)
\eeq

\noindent for causal impulse response functions $q_{ij}(t)$ and $p_{ik}(t)$ and vector of inputs $u(t) \in \RR^r$. Equation \eqref{yqp} defines a directed graph among observed variables and inputs with no self-loops and in which there is an edge from $y_j$ ($u_k$) to $y_i$ if and only if $q_{ij}(t) \not\equiv 0$ ($p_{ik}(t) \not\equiv 0$). The functions $q_{ij}(t)$ and $p_{ik}(t)$ therefore define both the topology and edge dynamics of the graph.

By taking the Laplace transform of \eqref{yqp}, the system is represented compactly in matrix form as follows:
\beq \label{YQP}
Y = QY + PU
\eeq

\noindent where $Y(s)$ and $U(s)$ are the Laplace transforms of $y(t)$ and $u(t)$ respectively; $Q(s)$ is a strictly-proper transfer matrix with entry $(i,j)$ equal to the Laplace transform of $q_{ij}(t)$ and diagonal entries equal to zero; $P(s)$ is a strictly-proper transfer matrix with entry $(i,k)$ equal to the Laplace transform of $p_{ik}(t)$. The couple $(Q,P)$ is termed the Dynamical Structure Function (DSF) and is uniquely defined for any partially-observed LTI system \cite{Gon08}.

A state-space realization of \eqref{yqp} can always be made in which the state vector is partitioned into \emph{manifest} states, corresponding to the observed variables $y_i$ and \emph{latent} states, which provide the dynamics of the impulse response functions $q_{ij}$ and $p_{ik}$. By observing more of the states, a different DSF will be obtained that represents the system in greater detail, hence the DSF may be regarded as a representation of the system at the \emph{resolution} ($\frac{p}{n}$, where $n$ is the total number of states) of the manifest states. Similarly, manifest states may be treated as latent in order to obtain a coarser DSF at a lower resolution.

% 1. Time-invariant system
% 2. Linearise about equilibrium
% 3. Discrete time using time-shift operator?
% 4. Subtract off w_{ii}s to remove self loops and preserve structure: q,p
% 5. Write in matrix form : Q,P
% 6. Mention this can be derived from partially-observed state space

% that interact via causal, Linear and Time-Invariant (LTI) functions. The states $y$ are driven by the vector of inputs $u(t) \in \RR^r$ also via LTI functions. This system defines a graph in which the nodes are the states $y$ and inputs $u$ and the edges are strictly proper transfer functions. Such a system has a state space realization $(A,B,C,D)$, where $D=0$ and $C = \bmat I & 0 \emat$ such that the state vector $x$ is partitioned into manifest ($y$) and latent states. Any system of this form also has a unique Dynamical Structure Function (DSF) representation: $(Q,P)$ defined in \cite{Gon08}. The matrices $Q \in \bbS^{p \times p}$ and $P \in \bbS^{p \times r}$ are strictly proper transfer matrices that describe both the edge topology of the network and open-loop dynamics of each of the edges, such that:
%\beq \label{YQPU}
%Y(s) = Q(s)Y(s) + P(s)U(s)
%\eeq
%
%\noindent where $Y(s)$ and $U(s)$ are the Laplace transforms of $y(t)$ and $u(t)$ respectively. The matrix $Q$ is constructed to have diagonal elements equal to zero to eliminate self loops. We pose the network reconstruction problem as obtaining $P$ and, in particular, $Q$ from $Y$ and $U$.

\subsection{Network Reconstruction}

We pose the network reconstruction problem as obtaining $(Q,P)$ from input/output data $(U,Y)$. Suppose $m$ independent experiments have been performed in each of which a different input (or set of inputs) has been applied. The largest number of independent experiments is $r$, the dimension of the inputs, so $m \leq r$. Denote the Laplace transforms of the inputs and outputs in the $i^{th}$ experiment as $U^{(i)}$ and $Y^{(i)}$ and concatenate these to form the following matrices:
\beqs
\bali
Y &:= \bmat Y^{(1)} & Y^{(2)} & \cdots & Y^{(m)} \emat \\
U &:= \bmat U^{(1)} & U^{(2)} & \cdots & U^{(m)} \emat
\eali
\eeqs

\noindent of dimension $p \times m$ and $r \times m$ respectively. Note that ${Y = QY + PU}$ and rearrange this to give:
\beq \label{YTQT}
\bmat Y^T & U^T \emat \bmat Q^T \\ P^T \emat  = Y^T
\eeq

\noindent where $\bmat Y^T & U^T \emat$ has dimension $m \times (p+r)$ and we wish to solve for $Q$ and $P$. By applying the vectorization operator we can write \eqref{YTQT} in the form of \eqref{Axb}:
\beq \label{AXB}
\bali
A \ &\leftarrow \ I \otimes \bmat Y^T & U^T \emat \\
x \ &\leftarrow \ vec\left( \bmat Q^T \\ P^T \emat \right) \\
b \ &\leftarrow \ vec\left(Y^T\right)
\eali
\eeq

\noindent where ${A(s) \in \CC^{M \times N}}$, $x(s) \in \CC^N$ and $b(s) \in \CC^M$ for $M = mp$ and $N = p(p+r)$. 

With no other information about the system, in order for \eqref{YTQT} to be well posed it is therefore required that:
\beq \label{MgN}
M\geq N \quad \Leftrightarrow \quad m \geq p+r \quad \Leftrightarrow \quad 0 \geq p
\eeq

\noindent since $r \geq m$, and hence additional information is always necessary. The assumption that rows of $Q$ and $P$ are sparse may be sufficient to ensure a unique sparse solution. In particular, if each row of $\bmat Q & P \emat$ is $k$-sparse then the solution to \eqref{YTQT} is unique if:
\beqs
M \geq 2kp \quad \Leftrightarrow \quad m \geq 2k
\eeqs

\noindent and the condition of Lemma \ref{lem:cs} is satisfied for $A$. However, this may not be the case due to the particular way in which the matrix $A$ is constructed in \eqref{AXB}.

Alternatively, we may assume some knowledge of how the inputs target the network, for example that $P$ is square ($r=p$) and diagonal. This defines an experimental setup in which each input is associated with a particular manifest state and affects it via the corresponding diagonal entry of $P$. By removing the $p-1$ known zero entries in each row of $P$ and the zero diagonal entries of $Q$, the condition for solution uniqueness \eqref{MgN} becomes:
\beqs
m \geq p \quad \Rightarrow \quad m=p
\eeqs

\noindent since $m \leq r = p$. This is the main result of \cite{Gon08} -- that with $P$ diagonal and no other \emph{a priori} information, $r=p$ is necessary and sufficient for solution uniqueness.

\section{Compressed Sensing with prior knowledge} \label{sec:prior}
Here we treat the problem of finding a solution $x \in \mathbb{R}^n$ to the following:
\begin{equation} \label{Axb2}
	Ax = b
\end{equation}

\noindent where $A \in \mathbb{R}^{m\times n}$ and $b \in \mathbb{R}^m$ are known, $m<n$ and where $x$ is partitioned into sparse and nonzero components. Without loss of generality, we can write \eqref{Axb2} as:
\beq \label{A12xb}
\bmat A_1 & A_2 \emat \bmat x_1 \\ x_2 \emat = b
\eeq

\noindent where $x_1 \in \RR^{n_1}$ satisfies $\|x_1\|_0 = k$, $x_2 \in \RR^{n_2}$ satisfies $\|x_2\|_0 = n_2$ and $n_1 + n_2 = n$. The vector $x$ is therefore ($k+n_2$)-sparse -- it satisfies $\|x\|_0 = k+n_2$.

This will be applied to the sparse reconstruction problem in the following section, where rows of $Q$ in \eqref{AXB} are sparse and part of $P$ is nonzero. Equation \eqref{Axb2} has a unique $(k+n_2)$-sparse solution by Lemma \ref{lem:cs} if $m \geq 2(k+n_2)$ and all subsets of $2(k+n_2)$ columns of $A$ are full rank. By making use of the known structure of $x$, we can solve for $x_1$ and $x_2$ separately and hence reduce the number of experiments, $m$, required.

This problem has been considered in the CS literature (see for example \cite{borries}, \cite{scarlett}) where it is well known that the number of experiments required for exact reconstruction can be reduced by $n_2$, the number of known nonzero entries of $x$. However, it is generally assumed that the matrix $A$ can be chosen to satisfy Lemma \ref{lem:cs}, which is not the case in our application. Here we derive a lemma analogous to Lemma \ref{lem:cs} that gives conditions for exact reconstruction.

\subsection{Conditions for Exact Reconstruction}
Take the $QR$ decomposition of $A_2 \in \RR^{m \times n_2}$:
\beq \label{A1QR}
A_2 = \bmat Q_1 & Q_2 \emat \bmat R_1 \\ 0 \emat
\eeq

\noindent where $\bmat Q_1 & Q_2 \emat \in \RR^{m \times m}$ is orthogonal and $R_1 \in \RR^{n_2 \times n_2}$ is upper triangular. Pre-multiply \eqref{A12xb} by $\bmat Q_1 & Q_2 \emat^T$:
\beq \label{QA12xb}
\bmat Q_1^TA_1  & R_1\\ Q_2^T A_1 & 0\emat \bmat x_1 \\ x_2 \emat = \bmat Q_1^Tb \\ Q_2^Tb \emat
\eeq

\noindent We can now solve first for $x_1$ using the second block row:
\beq \label{solvex2}
Q_2^TA_1x_1 = Q_2^Tb
\eeq

\noindent where $Q_2^TA_1 \in \RR^{(m-n_2) \times n_1}$ and $\|x_1\|_0 = k$. From Lemma \ref{lem:cs}, \eqref{solvex2} has a unique, $k$-sparse solution if $m-n_2\geq 2k$ and all subsets of $2k$ columns of $Q_2^TA_1$ are full rank. The number of experiments, $m$, required for solution uniqueness has therefore been reduced by the number of known components of $x$ from $2(k+n_2)$ to $2k +n_2$.

Given $x_1$, we may then solve for $x_2$ from the first block row of \eqref{QA12xb}:
\beqs
R_1x_2 = Q_1^T\left( b - A_1x_1 \right)
\eeqs

\noindent This has a unique solution if and only if $R_1$ is full rank, requiring $A_2$ to be full column rank. In this case, $x_2$ is given by:
\beq \label{solvex1}
x_2 = R_1^{-1}Q_1^T\left( b - A_1x_1 \right)
\eeq

\noindent The following lemma summarizes the above conditions for sparse solution uniqueness.
\begin{lem} \label{lem:cs2}
Suppose the sparsest solution to \eqref{A12xb} is known to have $\|x_1\|_0 = k$ and $\|x_2\|_0 = n_2$ and let the $QR$ decomposition of $A_2$ be given by \eqref{A1QR}. Then if $m \geq 2k + n_2$, if $A_2$ is full rank and if all subsets of $2k$ columns of $Q_2^TA_1$ are full rank, then the $(k+n_2)$-sparse solution to \eqref{A12xb} is unique.
\end{lem}

\begin{figure}[t]
\centering
\subfloat [] [Fixed $n_2$]{
\includegraphics[width=0.2\textwidth]{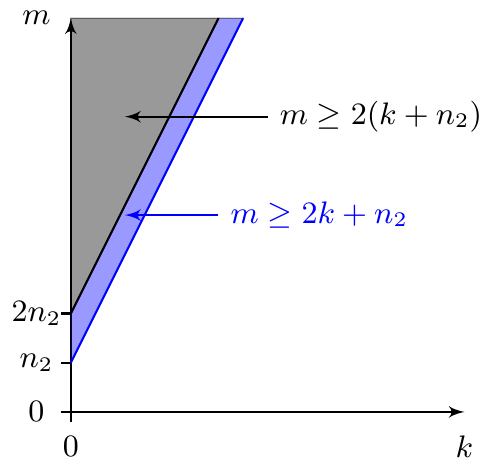}
}
\hspace{0.25cm}
\subfloat [] [Fixed $k$]{
\includegraphics[width=0.2\textwidth]{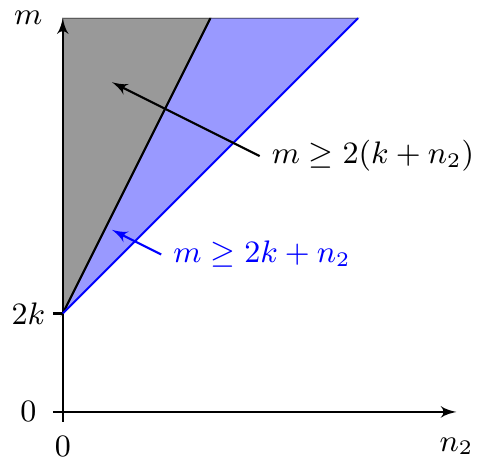}
}
\caption{Values of $m$ that satisfy Lemma \ref{lem:cs} (gray) and Lemma \ref{lem:cs2} (blue) for different values of $k$ (a) and $n_2$ (b). The blue region represents the increased number of problems that can be solved by incorporating prior information.}
\label{fig:mkn1}
\end{figure}

By incorporating prior knowledge the number of experiments $m$ required for $k$-sparse solution uniqueness is reduced by $n_2$, the number of known nonzero components. Equivalently, for a given number of experiments, the sparse solution for a higher value of $k$ may be unique. These ideas are illustrated in Fig. \ref{fig:mkn1}.

%% Discussion about using prior information
%\dph{Cut this?} It is also clear that uniqueness of a $(k+n_2)$-sparse solution to \eqref{Axb2} implies uniqueness of a $k$-sparse solution to \eqref{solvex2}, meaning that there is no detriment to using the knowledge that $\|x_2\|_0 = n_2$. To see this, suppose there are multiple $k$-sparse solutions $x_1^{(1)}$ and $x_1^{(2)}$ to \eqref{solvex2} with corresponding $x_2^{(1)}$ and $x_2^{(2)}$ found from \eqref{solvex1}. Since
%\beqs 
%x^{(1)} = \bmat x_1^{(1)} \\ x_2^{(1)} \emat \qquad \mathrm{and} \qquad x^{(2)} = \bmat x_1^{(2)} \\ x_2^{(2)} \emat
%\eeqs
%
%\noindent are both $(k+n_2)$-sparse and both solve \eqref{QA12xb}, they also solve \eqref{Axb2}. Hence non-uniqueness of sparse $x_1$ in \eqref{solvex2} implies non-uniqueness of sparse $x$ in \eqref{Axb2}.

%\dph{Alternative lemma:}
%\begin{lem}
%Suppose the sparsest solution to \eqref{A12xb} is known to have $\|x_1\|_0 = n_1$ and $\|x_2\|_0 = k$ and let $\phi_{2k}$ be a vector that indexes $2k$ columns of $A_2$. Then if $m \geq 2k + n_1$ and if $\bmat A_1 & A_2(:,\phi_{2k}) \emat$ is full rank for all $\phi_{2k}$, then the solution to \eqref{A12xb} is unique.
%\end{lem}

\section{Underdetermined Reconstruction} \label{sec:main}

Consider an unknown system defined by its Dynamical Structure Function (DSF): $(Q^0,P^0)$ where for inputs $U(s)$ and outputs $Y(s)$ we have:
\beqs
Y = Q^0Y + P^0U
\eeqs

\noindent Given some set of inputs and outputs we seek to identify $Q^0$ under the following assumptions.
\begin{ass} \label{ass:Pdiag}
The matrix $P^0$ is square, diagonal and full rank.
\end{ass}

\begin{ass} \label{ass:rlm}
The number of experiments $m$ is fewer than the number of measured states $p$.
\end{ass}

\begin{ass} \label{ass:Qsp}
The rows of $Q^0$ are $k$-sparse:
\beqs
\|Q^0(i,:)\|_0 \leq k < p \qquad \textrm{for} \qquad i=1, \ldots, p
\eeqs
\end{ass}

\noindent Assumption \ref{ass:Pdiag} asserts that the dimension of the input vector is equal to that of the manifest state vector ($r=p$) and that each input directly and uniquely affects one manifest state. Assumption \ref{ass:rlm} then states that we can only perform $m<p$ experiments, and hence without additional assumptions the reconstruction problem is ill posed.

By Assumption \ref{ass:Qsp}, each node in the graph defined by $Q^0$ has maximum in-degree of $k$ and the sparse solution to the reconstruction problem may be unique. In this case we seek sparse solutions $Q$ and diagonal $P$ to:
\beq \label{YTQT2}
\bmat Y^T & U^T \emat \bmat Q^T \\ P^T \emat  = Y^T
\eeq

\noindent This may be converted into the form of \eqref{A12xb} either by taking the vectorization operator as in \eqref{AXB} or by solving for each row of $\bmat Q & P \emat$ separately:
\beq \label{YTQTi}
\bmat Y^T & U^T(:,i) \emat \bmat Q^T(:,i) \\ P^T(i,i) \emat = Y^T(:,i)
\eeq

\noindent where $P(i,i) \neq 0$. The $(k+1)$-sparse solution to \eqref{YTQTi} is unique if the conditions of Lemma \ref{lem:cs2} are satisfied.

This problem is fundamentally different to typical compressed sensing applications in that we do not have free choice of the sensing matrix. First we will demonstrate that for a na\"{\i}ve choice of $U$, the condition of Lemma \ref{lem:cs2} may never be met; then we will consider the design of the sensing matrix by choice of $U$.

\subsection{Single Inputs}

Suppose that in each experiment only one input is applied, such that $U$ can be written as:
\beq \label{Udiag}
U = \bmat U_1 \\ 0 \emat
\eeq

\noindent where $U_1$ is square, diagonal and of dimension $m \times m$. Partitioning $Q$ and $P$ accordingly, the system equations can be written as:
\beq \label{pDSF}
\bmat Y_1 \\ Y_2 \emat =
\bmat Q_{11} & Q_{12} \\ Q_{21} & Q_{22} \emat
\bmat Y_1 \\ Y_2 \emat + 
\bmat P_{11} \\ 0 \emat U_1
\eeq

\noindent where $P_{11}$ is square and diagonal and the manifest states are therefore partitioned into \lq\lq perturbed\rq\rq \ ($Y_1$) and \lq\lq unperturbed\rq\rq \ ($Y_2$)\footnote{A state is unperturbed if it is not directly driven by an input -- it may still be indirectly affected via other perturbed states.}. Denote the DSF with respect to the $p$ manifest states $(Q,P)$ as the p-DSF; then by eliminating $Y_2$ from the right hand side of \eqref{pDSF} we can derive two further representations of the system.

From the first block row we can obtain the DSF with respect to only the $m$ perturbed states and denote this the m-DSF. This is a representation of the same system but at a lower resolution ($\frac{m}{n}$) of manifest states. Eliminating $Y_2$ gives:
\begin{equation}\label{eqmDSF}
\begin{aligned}
Y_1 &= \left( Q_{11} + Q_{12}\left( I - Q_{22} \right) ^{-1}Q_{21} \right)Y_1 + P_{11}U_1\\
 &=: \bar{Q}_{11}Y_1 + P_{11}U_1
\end{aligned}
\end{equation}

\noindent where in order to obtain a hollow $Q$ matrix (diagonal entries equal to zero) we must subtract the diagonal part of $\bar{Q}_{11}$ from both sides. Let ${D_{11} = \Diag(\bar{Q}_{11})}$, then:
\begin{equation}\label{eqmDSF2}
\begin{aligned}
Y_1 &= \bar{Q}_{11}Y_1 +D_{11}Y_1 - D_{11}Y_1+ P_{11}U_1\\
&= \left( I - D_{11} \right) ^{-1} \left( \left( \bar{Q}_{11} - D_{11} \right)Y_1 + P_{11}U_1 \right) \\
&=: \hat{Q}_{11}Y_1 + \hat{P}_{11}U_1
\end{aligned}
\end{equation}

\noindent The m-DSF is defined as $(\hat{Q}_{11},\hat{P}_{11})$, which treats only $Y_1$ as manifest states and $Y_2$ as additional latent states. This can be equivalently defined from a state-space realization by changing the partitioning of the state vector to reflect the change in manifest states. The m-DSF can therefore be obtained uniquely from $Y_1$ and $U_1$ since $U_{1}$ is square, diagonal and full rank \cite{Gon08}.

The second block row of \eqref{pDSF} gives:
\begin{equation}\label{eqPhi}
\begin{aligned}
Y_2 &= (I-Q_{22})^{-1}Q_{21}Y_1\\
&=: \hat{Q}_{21} Y_1
\end{aligned}
\end{equation}

\noindent where the transfer function $\hat{Q}_{21}$ describes causal relations from states in $Y_1$ to those in $Y_2$ that are direct in the sense that they do not involve other states in $Y_1$. The matrix $\hat{Q}_{21}$ is also identifiable from input/output data.
\begin{lem} \label{lem:ps}
A particular solution to the sparse network reconstruction problem \eqref{YTQT2} is:
\beqs
\hat{Q} = \bmat \hat{Q}_{11} & 0 \\ \hat{Q}_{21} & 0 \emat, \qquad \hat{P} = \bmat \hat{P}_{11} & 0 \\ 0 & 0 \emat
\eeqs

\noindent where $(\hat{Q}_{11},\hat{P}_{11})$ is defined in \eqref{eqmDSF2} and $\hat{Q}_{21}$ in \eqref{eqPhi}.
\end{lem}

\noindent Hence we can always construct at least one solution for $Q$, in which the unperturbed states, $Y_2$, have no outputs to any other measured states. Clearly if $\|\hat{Q}(i,:)\|_0 \leq k$ for any $i$ then the $k$-sparse solution to \eqref{YTQTi} is not unique.
\begin{rem}
Given $Q$ for which $\max_i \|Q(i,:)\|_0 = k$ and $\max_i \|\hat{Q}(i,:)\|_0 = \hat{k}$, it is possible that either $\hat{k} \leq k$ or $\hat{k} >k$. This result can be easily seen by example, and in the former case the $k$-sparse solution will not be unique.
\end{rem}

\noindent The assumption of sparsity for any network must therefore be firmly justified \emph{a priori} for the representation in question -- when viewed at a different resolution of manifest states, the sparsity of the network can change.

%The system in Example \ref{egring} in the following section satisfies $\|Q(i,:)\|_0 = \|\hat{Q}(i,:)\|_0$ for all $i$.

\subsection{Constraints on $Q$}

Assume single inputs of the form \eqref{Udiag} have been applied and partition $Q$ and $\hat{Q}$ of Lemma \ref{lem:ps} as follows:
\beqs
Q = \bmat Q_1 & Q_2 \emat, \qquad \hat{Q} = \bmat \hat{Q}_1 & 0 \emat, \qquad \hat{Q}_1 := \bmat \hat{Q}_{11} \\ \hat{Q}_{21} \emat
\eeqs

\noindent such that $Q_1$ has the same dimension as $\hat{Q}_1$. From the identifiable quantity $\hat{Q}_1$ we may infer something about whether entries of $Q$ are zero or not.
\begin{lem} \label{lem:Qcon}
For every $i \neq j$, if $\hat{Q}_1(i,j) \neq 0$, then
\beqs
Q_1(i,j) \neq 0 \qquad \mathrm{or} \qquad Q_2(i,k)\hat{Q}_{21}(k,j) \neq 0
\eeqs

\noindent for some $k$. Else $\hat{Q}_1(i,j) = 0$, then
\beqs
Q_1(i,j) = 0 \aaand Q_2(i,k)\hat{Q}_{21}(k,j) = 0
\eeqs

\noindent for every $k$ unless the graph defined by $Q$ contains multiple paths that sum to zero.
\end{lem}

\begin{pf}
From the definition of $\hat{Q}_{11}$ in \eqref{eqmDSF2}, for $i\neq j$:
\beqs
\hat{Q}_{11}(i,j) \neq 0 \iiff \bar{Q}_{11}(i,j) \neq 0
\eeqs

\noindent Then from \eqref{eqmDSF2} and \eqref{eqPhi} we have:
\beq \label{Q1212}
\bmat \bar{Q}_{11} \\ \hat{Q}_{21} \emat =
\bmat Q_{11} \\ Q_{21} \emat + \bmat Q_{12} \\ Q_{22} \emat \hat{Q}_{21}\\
= Q_1 + Q_2 \hat{Q}_{21}
\eeq

\noindent If $\hat{Q}_1(i,j) \neq 0$, the result now follows directly from \eqref{Q1212}. If $\hat{Q}_1(i,j) = 0$ it is possible that $Q_1 + Q_2 \hat{Q}_{21} = 0$ and the condition of the lemma not be satisfied. However, $\hat{Q}_1(i,j)$ comprises the sum of a direct link $Q_1(i,j)$ and a sum of paths:
\beqs
Q_2(i,:)\hat{Q}_{21}(:,j) = Q_2(i,:)\left(I-Q_{22}\right)^{-1}Q_{21}(:,j)
\eeqs

\noindent in $Q$. Hence this case necessitates that this direct link and all the paths sum to zero.
\end{pf}

The case of multiple paths summing to zero in $Q$ is considered unlikely to occur in practice. Lemma \ref{lem:Qcon} may be interpreted as low resolution structure $\hat{Q}_{1}$ having to be consistent with the structure of $Q$ -- a connection that exists in $\hat{Q}_1$ must be preserved (directly or indirectly) in $Q$. The following example illustrates how this can be used to place constraints on the unknown $Q$.

\begin{figure}[t]
\centering
\subfloat[][$Q$]{
\centering
\includegraphics{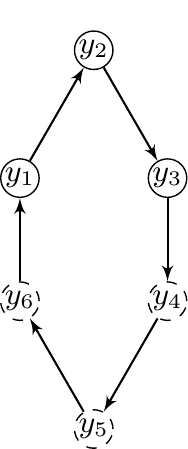}
}
$\hspace{2cm}$
\subfloat[][$\hat{Q}$]{
\centering
\includegraphics{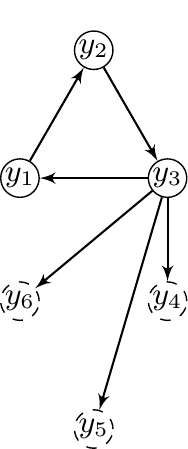}
}
\caption{Network with $p=6$ and $m=3$ in which the first three manifest states are perturbed. Solid circles denote perturbed states ($Y_1$), dashed circles unperturbed states ($Y_2$) and arrows denote nonzero entries of (a) $Q$ and (b) $\hat{Q}$.}
\label{fig:eg1}
\end{figure}

\begin{myeg} \label{egring}
Consider the network of Fig. \ref{fig:eg1}(a) with $p=6$ and $m=3$. The matrix $Q$ and the particular solution $\hat{Q}$ from Lemma \ref{lem:ps} are:
\beqs
Q = \bmat
0 & 0 & 0 & 0 & 0 & \times \\
\times & 0 & 0 & 0 & 0 & 0 \\
0 & \times & 0 & 0 & 0 & 0 \\
0 & 0 & \times & 0 & 0 & 0 \\
0 & 0 & 0 & \times & 0 & 0 \\
0 & 0 & 0 & 0 & \times & 0 \emat
, \qquad
\hat{Q} = \bmat
0 & 0 & \times & 0 & 0 & 0 \\
\times & 0 & 0 & 0 & 0 & 0 \\
0 & \times & 0 & 0 & 0 & 0 \\
0 & 0 & \times & 0 & 0 & 0 \\
0 & 0 & \times & 0 & 0 & 0 \\
0 & 0 & \times & 0 & 0 & 0 \emat
\eeqs

\noindent where $\times$ denotes a nonzero entry. The matrix $\hat{Q}$ is a valid $1$-sparse solution to \eqref{YTQTi} and is shown in Fig. \ref{fig:eg1}(b). Using Lemma \ref{lem:Qcon} we can obtain a matrix $Q^c$ which describes constraints on the solution set of $Q$:
\beqs
Q^c = \bmat
0 & 0 & ? & ? & ? & ? \\
\times & 0 & 0 & 0 & 0 & 0 \\
0 & \times & 0 & ? & ? & ? \\
0 & 0 & ? & 0 & ? & ? \\
0 & 0 & ? & ? & 0 & ? \\
0 & 0 & ? & ? & ? & 0 \emat
\eeqs

\noindent where $?$ denotes an unknown entry. Approximately half of the structure of $Q$ can therefore be found, but the $1$-sparse solution is unique for rows two and three only. For any such ring network with $m<p$, only rows $2,\ldots,m$ have a unique sparse solution and hence the solution for $Q$ is unique if and only if $m=p$.
\end{myeg}

\noindent With single inputs described by \eqref{Udiag}, the assumption of sparsity alone is not sufficient to ensure uniqueness of the sparse solution $(Q,P)$ to the network reconstruction problem. This is illustrated by Example \ref{egring} in which $Q$ is $1$-sparse (unit in-degree) and has multiple $1$-sparse solutions for any $m<p$. Applying inputs in this manner does however yield lower resolution structural information that can impose constraints on $Q$ by Lemma \ref{lem:Qcon}. In addition, we have not considered the minimality or stability of the solutions -- for small problems it may be clear, for example, that the solution of minimal dimension is unique. Next we consider the design of $U$ to satisfy Lemma \ref{lem:cs2} and hence ensure sparse solution uniqueness.

%%%%%

%%%%%

\subsection{Experiment Design} \label{subs:expt}
\begin{figure}[t]
\centering
\subfloat[][Experiment 1]{
\includegraphics{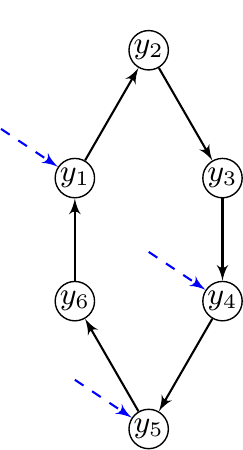}
}
\hspace{0.1cm}
\subfloat[][Experiment 2]{
\includegraphics{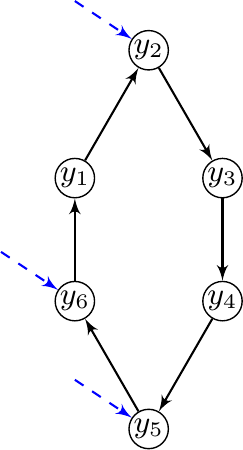}
}
\hspace{0.1cm}
\subfloat[][Experiment 3]{
\includegraphics{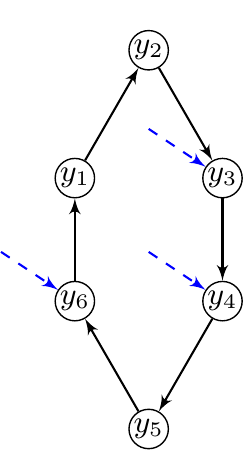}
}
\caption{The network of Fig. \ref{fig:eg1}(a) with $p=6$, $m=3$ and three inputs applied in each experiment denoted by the blue dashed arrows.}
\label{fig:eg2}
\end{figure}

The problem with the diagonal inputs of \eqref{Udiag} is that the unperturbed states have no variation independent of their parent states. As a result, subsets of columns of $Y^T$ in \eqref{YTQTi} may not be full rank and hence may not satisfy Lemma \ref{lem:cs2}. Figure \ref{fig:eg2} shows the same network as in Fig. \ref{fig:eg1}(a) with three inputs applied in each of the three experiments. It is straightforward to verify that any system with this structure generically satisfies Lemma \ref{lem:cs2} for $k=1$ for any particular $s=j\omega$ and therefore has a unique $1$-sparse solution. However, without first knowing the structure, one would not be able to design such inputs.

%In the above example, three perturbations per experiment suffice, but it is also important which perturbations are applied in which experiment. In general when designing experiments we must make two choices: how many perturbations and which states to perturb in each experiment.

Perturbation design was considered in \cite{tegner} for fully-observed linear systems using steady-state data. An algorithm is presented in which perturbations are first applied at random until every state has been perturbed in at least one experiment. Then all solutions of a certain sparsity consistent with the data are constructed and a measure of variance for each state is obtained based on how much the outputs of this state vary across all the solutions. The state(s) with the highest variance are then perturbed and the procedure iterated. Simulation results suggested that the more states that are perturbed in each experiment, the fewer experiments are needed.

%The number of experiments needed for a network of $n$ states with maximum in-degree of $k$ was found empirically to be $\mathcal{O}(k\log n)$.
%In practice inputs may be expensive to apply and if the inputs cannot be measured this may make the solution less sparse, as discussed in Section \ref{sec:inputs}. 

%%%%%
%%%%%

A similar procedure is presented in \cite{karp} for acyclic Boolean networks. At each stage of the algorithm, the next experiment is selected as that which maximizes the decrease in entropy in terms of reducing the set of possible solutions. The set of possible experiments to choose from is taken as given. In \cite{Gon08} and \cite{sontag} targeted inputs are considered but the number of experiments must be equal to the number of manifest states. In the former, single inputs are considered such that each state is perturbed in turn; in the latter the opposite: every state except one must be perturbed in each experiment.

%, and again the choice of how many states should be perturbed in each experiment is not treated.

Three strategies are presented here along similar lines to \cite{tegner}; in each, experiments are performed iteratively until the solution at a given level of sparsity is unique, according to Lemma \ref{lem:cs2}. At each iteration, a fixed number of inputs, $l$, are chosen according to:
\begin{enumerate}
\item Random -- inputs are chosen with equal probability
\item Biased Random -- inputs are chosen at random with a bias towards those that have been applied the least in previous experiments
\item Targeted -- inputs are chosen to target rank-deficient subsets of columns of the sensing matrix in Lemma \ref{lem:cs2}. Each column of the sensing matrix corresponds to a manifest state and therefore to an input; having identified a deficient subset, the input in this subset that has been applied the least in previous experiments is selected
\end{enumerate}

\noindent The third approach will be seen in Section \ref{sec:sim} to require the fewest experiments on average to ensure a unique solution; it does however incur the additional computational cost of searching for rank-deficient subsets.

\section{Numerical Simulations} \label{sec:sim}

\subsection{Experimental Design for Solution Uniqueness} \label{subs:simex}

Here we compare in simulation the number of experiments required for a unique solution by the three strategies of Section \ref{subs:expt}. Random networks of $p=20$ measured states were generated with maximum in-degree sparsity of $k=1,2,3,4$. For each network, we performed experiments with exactly $l$ step inputs applied, for $l = 1,\ldots,p$ following each of the three strategies. In each case, the steady-state response of the network was used to assess whether the conditions of Lemma \ref{lem:cs2} were satisfied. The results are given here for $k=2$; the same trends were observed for the other values of $k$.

Figure \ref{fig:sim1} shows the average number of experiments needed for Lemma \ref{lem:cs2} to be satisfied over 100 trials for $k=2$. If only one input is applied in each experiment, the maximum number of experiments is always required and CS hence offers no improvement. Substantial reduction in the number of experiments is observed if more than one input can be applied, and applying these in a biased or targeted manner is more effective than applying them at random. The best results are seen for the targeted approach, although as mentioned this incurs higher computational cost. For a sufficiently large number of inputs per experiment, the number of experiments required for a unique solution is determined by $m = 2k +1$, independent of the strategy.

\subsection{Exact Reconstruction using Basis Pursuit}
\begin{figure}[t]
\centering
\includegraphics{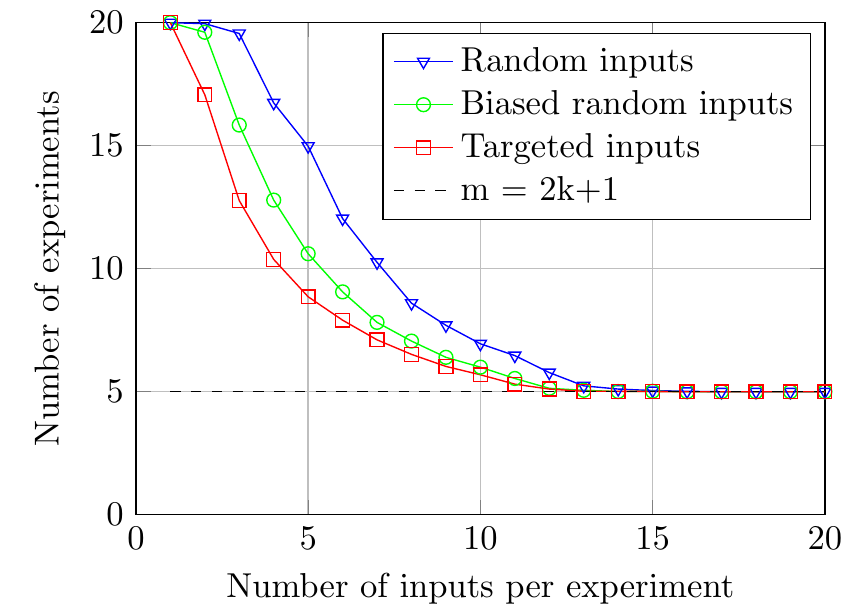}
\caption {Mean number of experiments needed for solution uniqueness (Lemma \ref{lem:cs2}) for different numbers of inputs per experiment for three different strategies. The mean was taken over 100 random networks of $p=20$ manifest states with maximum in-degree sparsity of $k=2$.}
\label{fig:sim1}
\end{figure}

For 100 random networks of the type treated in Fig. \ref{fig:sim1} with four inputs per experiment, we attempt to solve for the 2-sparse solution using basis pursuit ($l_1$ minimization). A sufficient condition for the success of basis pursuit is that the coherence of the sensing matrix is sufficiently small. In \cite{sanad}, lower bounds on the network coherence of simple networks were derived and also observed in simulation as the number of measurements was increased. In particular the lower bound was increasing in the magnitude of the parameters of the impulse response functions. A similar phenomenon was observed here: increasing the steady-state gain of the entries of $Q$ and $P$ increased the coherence.

\begin{figure}[t]
\centering
\includegraphics{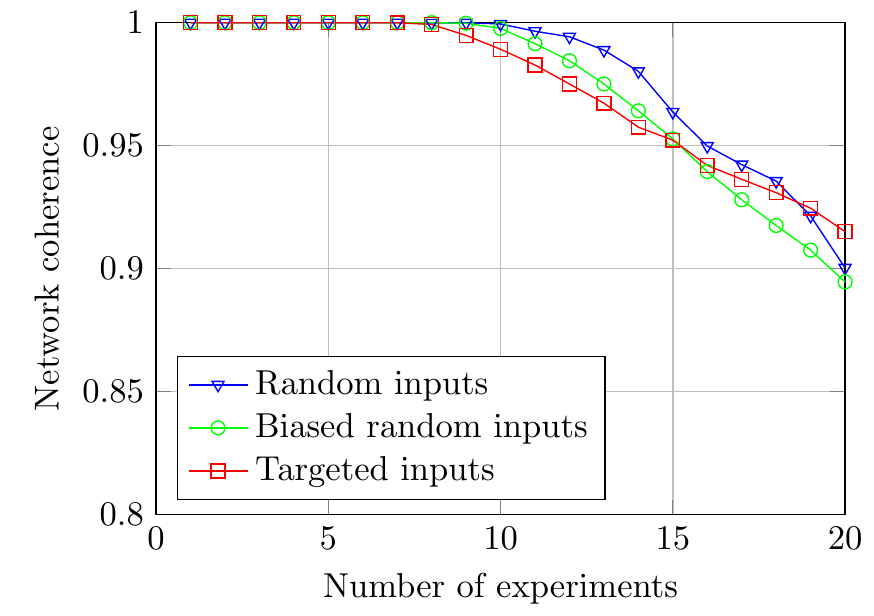}
\caption{Mean coherence of the sensing matrix for 100 random networks of $p=20$ manifest states with maximum in-degree sparsity of $k=2$ and four inputs applied in each experiment.}
\label{fig:BPcoh}
\end{figure}

\begin{figure}[t!]
\centering
\includegraphics{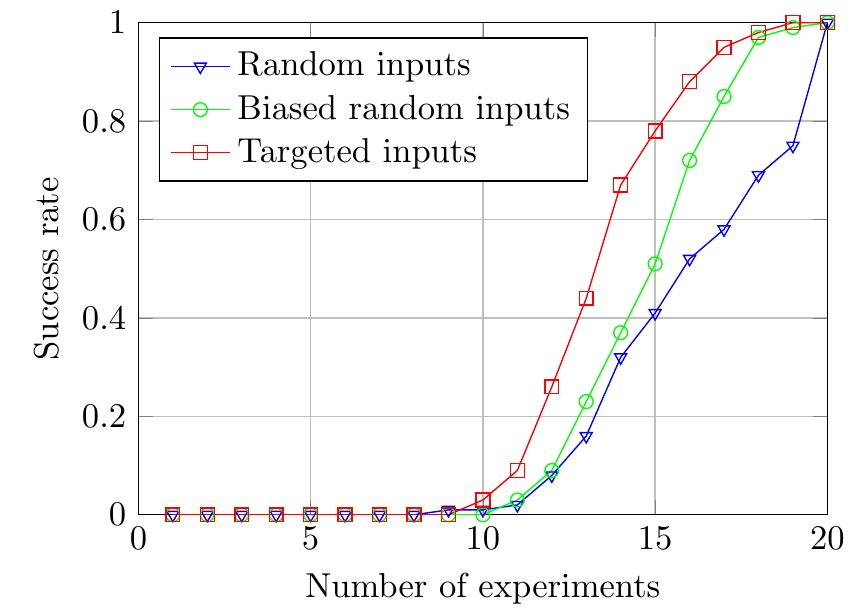}
\caption{Mean success rate using basis pursuit for the 100 networks of Fig.~\ref{fig:BPcoh}.}
\label{fig:BPsuc}
\end{figure}

For $|Q_{ij}(0)|, |P_{ii}(0)| < 0.5$, Fig. \ref{fig:BPcoh} shows the mean coherence of the sensing matrix against the number of experiments. Fig. \ref{fig:BPsuc} then shows the success rate of each of the three strategies at recovering the entire network from steady-state data using basis pursuit. The success rate is one if all links are recovered with no false positives and is zero otherwise; for this quite stringent metric, the performance is promising, particularly for the targeted experiments. The relative performance of each of the strategies is consistent with the solution uniqueness results of Fig. \ref{fig:sim1}. The actual number of experiment required for exact reconstruction using basis pursuit can be seen to be higher than that required for a unique solution, as is normally the case \cite{Candes2007}.

%%%%%

\section{Conclusions} \label{sec:conc}
We have investigated Compressed Sensing (CS) as a tool for reconstructing dynamical networks from data with a deficient number of experiments. The application differs from the typical use of CS in that the sensing matrix cannot be freely chosen and exact reconstruction hence necessitates appropriate choice of experiments. We provide a formulation of the problem that incorporates prior knowledge and present strategies for experimental design to attain sparse solution uniqueness with fewer experiments. 

The problem is motivated by biological applications where data are typically scarce; in this context we provide an algorithm for reconstruction at particular frequencies, for example at steady state, where a small number of data points suffice. It is also straightforward to apply this approach when the inputs are unknown, as may be the case in practice. Simulations demonstrate that the problem may be solved efficiently using a convex relaxation, such as basis pursuit, if the network coherence metric is sufficiently small. We also observe previously identified lower bounds on the network coherence, which highlight a potential limitation of the use of CS for this application.

\begin{ack}
This research was supported by the Engineering and Physical Sciences Research Council under Grant EP/G066477/1 and by the NIH NCI under the ICBP and PS-OC programs (5U54CA112970-08).
\end{ack}

\bibliographystyle{unsrt}
\bibliography{./mycsbib}
\end{document}